\documentclass[journal,10pt,onecolumn,draftclsnofoot,]{IEEEtran}

\usepackage[T1]{fontenc}

\ifCLASSINFOpdf
\else
\fi
\usepackage{amsmath}
\usepackage[cmintegrals]{newtxmath}

\usepackage{epsfig}
\usepackage{color}

\usepackage{amsmath} 
\usepackage[normalem]{ulem}
\newtheorem{definition}{Definition}[section]
\newtheorem{lemma}[definition]{Lemma}
\newtheorem{theorem}[definition]{Theorem}

\newtheorem{corollary}[definition]{Corollary}
\newtheorem{remarkth}[definition]{Remark}
\newtheorem{exampleth}[definition]{Example}
\newenvironment{remark}{\begin{remarkth}\upshape}{\hfill$\diamond$\end{remarkth}}
\newenvironment{example}{\begin{exampleth}\upshape}{\hfill$\diamond$\end{exampleth}}
\renewcommand{\emph}[1]{{\bfseries\itshape{#1}}}

\newcommand{\R}{\mathbb{R}}      

\makeatletter

\begin{document}
%
\title{Dynamic interpolation for  obstacle avoidance on Riemannian manifolds}
%
%
%

\author{Anthony Bloch,
        Margarida Camarinha, 
        and~Leonardo Colombo
\thanks{ A. Bloch is with Department of Mathematics, University of Michigan, 530 Church St. Ann Arbor, 48109, Michigan, USA.
        {\tt\small abloch@umich.edu}}
\thanks{M. Camarinha is with CMUC -- Centre for Mathematics of the University of Coimbra, Department of Mathematics, University of Coimbra
 3001-501 Coimbra, Portugal.
	    {\tt\small mmlsc@mat.uc.pt} }
\thanks{L. Colombo is with Instituto de Ciencias Matem\'aticas, Consejo Superior de Investigaciones Cient\'ificas. Calle Nicol\'as Cabrera 13-15, Campus UAM, Cantoblanco, 28049, Madrid, Spain. {\tt\small leo.colombo@icmat.es}}
}

\maketitle

\begin{abstract}
This work is devoted to studying dynamic interpolation for obstacle avoidance. This is a problem that consists of minimizing a suitable energy functional among a set of admissible curves subject to some interpolation conditions. The given energy functional depends on velocity, covariant acceleration and on artificial potential functions used for avoiding obstacles.

We derive first order necessary conditions for optimality in the proposed problem; that is, given interpolation and boundary conditions we find the set of differential equations describing the evolution of a curve that satisfies the prescribed boundary values, interpolates the given points and is an extremal for the energy functional.

We study the problem in different settings including a general one on a Riemannian manifold and a more specific one on a Lie group endowed with a left-invariant metric. We also consider a sub-Riemannian problem. We illustrate the results with examples of rigid bodies, both planar and spatial, and underactuated vehicles including a unicycle and an underactuated unmanned vehicle.\end{abstract}


%
\IEEEpeerreviewmaketitle

\section{Introduction}

Motion planning is an important task in numerous engineering fields such as air traffic control, aeronautics, robotics and computational anatomy. In the last few decades research in the calculus of variations and optimization has provided several methods for trajectory planning with  smooth interpolation by means of the study of higher-order  variational problems \cite{Altafini}, \cite{CMdD}, \cite{CroSil:91}, \cite{GHR10}, \cite{BlochHussein}, \cite{MacSil:2010},  \cite{Noa:89}, \cite{Zefran}.

In this class of problems the aim is to plan smooth  trajectories  passing through determined points at specific times. To achieve this goal, second order variational methods have been successfully used, providing interpolating curves, the so-called Riemannian  cubic splines and cubic splines in tension. These curves are better interpolating curves than geodesics, which usually fail  smoothness requirements for the trajectories. The variational problems consist of minimizing an energy functional, depending on the  covariant acceleration and velocity,  among a set of admissible curves interpolating a given set of points. Dynamic interpolation problems were initially studied by Crouch and Jackson \cite{Jackson} for applications in aeronautics and further explored from the geometric point of view by Noakes et al. \cite{Noa:89} and Crouch and Silva Leite \cite{CroSil:91}, \cite{CroSil:95}, \cite{SCC00}.

Crouch and Silva Leite \cite{CroSil:95} started the study of geometric properties of cubic polynomials on Riemannian manifolds, in particular on compact and connected Lie groups and symmetric spaces. Further extensions were developed by Bloch and Crouch \cite{BlCr}, \cite{blochcrouch} in the context of sub-Riemannian geometry. Sub-Riemannian problems are variational in nature with additional nonholonomic constraints. Although  nonholonomic variational problems do not give the correct approach to the study of nonholonomic mechanical systems, they provide an appropriate formulation for optimal control of kinematic underactuated control systems, including  rigid body control systems, which are of interest in the areas of aeronautics and robotics, since they are kinematic models for aerospace and underwater unmanned vehicles \cite{Biggs}, \cite{Noami2}, \cite{Vershik}.

Over the last few decades many authors have studied the problem of trajectory planning of autonomous vehicles in the presence of static obstacles in the workspace. Artificial potential functions \cite{Khatib1986} (as for instance, a Coulomb potential) have frequently been used for avoiding collision with obstacles, playing a fundamental role in these studies. These functions are created to simulate a fictitious repulsion from determined obstacles given by regions in the configuration space. This approach has been studied by Khabit for robotic manipulators (see \cite{Khatib1986} and references therein), and further studied by Koditschek \cite{K88} in the context of mechanical systems and Fiorelli and Leonard \cite{Leonard2} for multi-agent formation. The mathematical foundations for the existence of such a smooth functions on any smooth manifold can be found in the works of Smale \cite{Smale}, \cite{HS}.


In this paper, we aim to generate trajectories interpolating prescribed points and avoiding multiple obstacles in the workspace via the study of a second order variational  problem on a Riemannian manifold $M$ by an extension of the results presented in \cite{BlCaCoCDC} for variational obstacle avoidance without interpolation points. We call this problem \textit{dynamic interpolation for obstacle avoidance}. We study the problem in different scenarios: a general one on a Riemannian manifold and a more specific one on a Lie group endowed with a left-invariant metric, which is the appropriate setting for the examples we are interested in. We also consider the corresponding sub-Riemannian problem where we must deal with constraints on the velocities defined by a non-integrable distribution on $M$. We illustrate  the results with the examples of rigid bodies, both planar and spatial, and underactuated vehicles including the unicycle.

Our design for interpolation among piecewise smooth trajectories is independent of the knowledge of mass and inertia coefficients, leading to robustness for parameter uncertainty. Moreover, the solution for the sub-Riemannian problem allows for vehicle designs that include fewer actuators than is typical, leading to lighter, less costly design.

The framework proposed here for dynamic interpolation with obstacle avoidance of kinematic control systems on Lie groups endowed with a left-invariant metric should be useful for control design for a general
class of systems including spacecraft and underactuated vehicles.  In  general,  the  configuration space  for  these  systems  is  globally  described  by  a  matrix Lie  group  making this  model a  natural choice for  the  controlled system.  The  Lie  group description  gives  rise  to  coordinate-free expressions for the dynamics determining the behavior of the  system.  When  systems  on  Lie  groups  are  left  invariant, there  is  a  natural ``globalization''  of  solutions, that  is,  even  if  we exploit  local  charts  to  design local  maneuvers,  the Lie group formalism  allows us  to move over the entire configuration space without reformulating the controls. This is because we can
always left translate back to the identity of the group.

We extend previous results on dynamic interpolation on Riemannian manifolds to include obstacle avoidance. The main advantage of this approach is that it can be used to design global motions for many systems of practical interest where obstacle avoidance is necessary. The results of this work can be applied to a wide range of problems in systems and control such as spacecraft docking \cite{HB2018}, quantum control \cite{deAlessandro}, control of quadrotor UAVs \cite{lee}, multi-agent systems \cite{liu}, DNA structures \cite{Goyal}, control of marine cables/rods \cite{Bretl} and constrained under-actuated spacecraft \cite{Biggs17}, \cite{Biggs16}, \cite{Hu} among others.

The structure of the paper is as follows.

In section 2 we review the main topics  we use from Riemannian geometry and consider the variational  obstacle avoidance problem on Riemannian manifolds.

In section $3$ we introduce interpolation points into the previous framework in order to formulate the dynamic interpolation for obstacle avoidance problem. We
derive first order necessary conditions for optimality. By introducing the structure of a left-invariant Riemannian metric we study the problem on a Lie group. In section $4$ we extend our analysis to the sub-Riemannian situation characterizing stationary paths for an extended action integral with  constraints in velocities.

We apply the results to the dynamic interpolation for obstacle avoidance problems of several rigid bodies type systems, both planar and spatial, on the Lie groups $SE(2)$  and $SE(3)$, respectively, in section $3$, and underactuated vehicles  in section $4$. Final comments and ongoing work are discussed at the end of the paper.


\section{Variational obstacle avoidance problem on a Riemannian manifold}

\subsection{Preliminaries on Riemannian Geometry}

Let $M$ be a smooth ($\mathcal{C}^{\infty}$) \textit{Riemannian
manifold} with the \textit{Riemannian metric} denoted by
$\langle\cdot,\cdot\rangle$.  The length of a tangent vector is determined by its norm,
$||v_x||=\langle v_x,v_x\rangle^{1/2}$ with $v_x\in T_xM$, for each point $x\in M$.

A \textit{Riemannian connection} $\nabla$ on $M$, is a map that assigns to any two smooth vector fields $X$ and $Y$ on $M$ a new vector field, $\nabla_{X}Y$, called the \textit{covariant derivative of
$Y$ with respect to $X$}. For the properties of $\nabla$, we refer the reader to \cite{Boothby, bookBullo,Milnor}.

Consider a vector field $W$  along a curve $x$ on $M$. The $s$th-order covariant  derivative along $x$  of $W$ is denoted by $\displaystyle{\frac{D^{s}W}{dt^{s}}}$, $s\geq 1$. We also denote by $\displaystyle{\frac{D^{s+1}x}{dt^{s+1}}}$ the $s$th-order covariant derivative along $x$ of the velocity vector field of $x$, $s\geq 1$.

Given vector fields $X$, $Y$
and $Z$ on $M$, the vector field $R(X,Y)Z$ given by \begin{equation}\label{eq:CurvatureTensorDefinition}
R(X,Y)Z=\nabla_{X}\nabla_{Y}Z-\nabla_{Y}\nabla_{X}Z-\nabla_{[X,Y]}Z
\end{equation}  defines the \textit{curvature tensor} of $M$, where $[X,Y]$ denotes the \textit{Lie bracket} of the vector fields $X$ and $Y$. $R$ is trilinear in $X$, $Y$ and $Z$ and a tensor of type $(1,3)$. Hence for vector fields $X,Y,Z,W$ on $M$ the curvature tensor satisfies (\cite{Milnor}, p. 53)
\begin{equation}\label{curvformula}\langle R(X,Y)Z,W\rangle=\langle R(W,Z)Y,X\rangle.\end{equation}

\begin{lemma}[\cite{Boothby}, \cite{BlCr}]\label{lemma2}
Let $\omega$ be a one form on $(M,\langle\cdot,\cdot\rangle)$. The \textit{exterior derivative} of a one form $\omega$ is given by $$d\omega(X,Y)=X\omega(Y)-Y\omega(X)-\omega([X,Y])$$ for all vector fields $X,Y$ on $M$.
In particular, if $\omega(X)=\langle W,X\rangle$ it follows that
\begin{equation}\label{difomega}
d\omega(X,Y)=\langle\nabla_{X}W,Y\rangle-\langle\nabla_{Y}W,X\rangle.
\end{equation}
\end{lemma}

Let $\Omega$ be the set of all ${\mathcal C}^1$ piecewise smooth curves $x:[0, T]\to M$ in
$M$ such that $x(0)$, $x(T)$, $\frac{dx}{dt}(0)$ and $\frac{dx}{dt}(T)$ are fixed. The set $\Omega$ is called the
\textit{admissible set}.

For the class of curves in $\Omega$, we introduce the ${\mathcal C}^1$ piecewise smooth \textit{one-parameter admissible variation} of a curve $x\in\Omega$ by $
\alpha : (-\epsilon , \epsilon ) \times [0,T]  \rightarrow  M ;(r,t) \mapsto \alpha (r,t)=\alpha_r(t)$
that verify $\alpha_0(t)=x(t)$ and $\alpha_r\in\Omega$, for each $r\in (-\epsilon , \epsilon )$.

The \textit{variational vector field}  associated to an admissible variation $\alpha$ is a ${\mathcal C}^{1}$-piecewise smooth vector field along $x$ defined by $$X(t)=\frac{D}{\partial r}\Big{|}_{r=0}\alpha(r,t)\in T_{x}\Omega,$$
verifying  the boundary conditions
\begin{equation}\label{3.6}
X(0)=0,\quad   X(T)=0,\quad
\frac{DX}{dt}(0)=0,\quad   \frac{DX}{dt}(T)=0,\end{equation} where the tangent space  of  $\Omega$ at $x$ is the
 vector space $T_{x}\Omega$  of all ${\mathcal C}^{1}$ piecewise smooth vector
fields $X$ along $x$  verifying the boundary conditions (\ref{3.6}).

\begin{lemma}[\cite{Milnor}, p.$52$]\label{lemma} The one-parameter variation satisfies
$$\frac{D}{\partial r}\frac{D^2\alpha}{\partial
t^2}=\frac{D^2}{dt^2}\frac{\partial \alpha}{\partial r}+R\left(\frac{\partial \alpha}{\partial r}
,\frac{\partial \alpha}{\partial t}\right)\frac{\partial
\alpha}{\partial t}$$ where $R$ is the curvature tensor.
\end{lemma}

\subsection{The variational obstacle avoidance problem}

Let $T$, $\sigma$ and $\tau$  be positive real numbers, $(p_0,v_0)$, $(p_T,v_T)$ points in $TM$  and $S$ a regular submanifold of $M$. Consider the set $\Omega$ of all ${\mathcal C}^{1}$ piecewise smooth curves  on $M$, $x:[0,T]\rightarrow M$, verifying the boundary conditions \begin{equation}\label{3.1}
x(0)=p_0, \quad  x(T)=p_T,\quad \frac{dx}{dt}(0)=v_0, \quad \frac{dx}{dt}(T)=v_T,\end{equation} and define the functional $J$ on $\Omega$ given by
\begin{equation}\label{3.2}
J(x)=\int_{0}^{T}\frac{1}{2}\left(\Big{\|}\frac{D^2x}{dt^2}(t)\Big{\|}^2+
 \sigma \Big{\|}\frac{ dx}{dt}(t)\Big{\|}^2+V(x(t))\right)dt.
\end{equation}

\noindent This functional is  given by a weighted combination of the velocity and  covariant acceleration of the curve $x$ regulated by the
parameter $\sigma$, together with an artificial potential function $V:M\to\mathbb{R}$ used to avoid the obstacle. The obstacle is described by a region in $M$ bounded by $S$.

The potential function $V$ is an artificial smooth (or at least $C^{2}$) potential function
associated with a fictitious force inducing a repulsion from $S$. We consider $S$ to be the regular zero level set defined by a scalar valued smooth function  $f$, for instance, used to describe obstacles as circles in the plane for 2D vehicles, spheres or ellipsoids in the space for 3D vehicles, and orientations in the space for 3D rigid bodies.

To avoid collision with obstacles we introduce a potential function $V$ defined as the inverse value of the function $f$. The function $V$  goes to infinity  near the obstacle and decays to zero at some positive level set far away from the obstacle $S$.
This ensures that  such an optimal trajectory does not intersect $S$.  The use of artificial potential functions to avoid collision was introduced by Khatib (see \cite{Khatib1986} and references therein) and further studied by Koditschek \cite{K88}.

 \textbf{Problem 1}: The \textit{variational obstacle avoidance problem} consists of minimizing the functional $J$ on $\Omega$.

In order to extremize the functional $J$ on  the set $\Omega$ one needs to  compare the value of $J$ at a curve $x\in \Omega$ to the value of $J$ at a
nearby curve $\tilde{x}\in \Omega$, using one-parameter admissible
variations of $x$ in $\Omega$. We recently proved in \cite{BlCaCoCDC} the following result.
\begin{theorem}\label{CDC}
If  $x \in \Omega$ is an extremizer of $J$, then $x$ is smooth on $[0,T]$ and verifies
\begin{equation} \label{3.8}
\frac{D^{4}x}{dt^{4}}+R\left(\frac{D^2x}{dt^2},\frac{dx}{dt}\right)\frac{dx}{dt}- \sigma
\frac{D^2x}{dt^2}+ \frac{1}{2}\mbox{grad }V(x) = 0.
\end{equation}
\end{theorem}

\begin{remark}
In the absence of obstacles, we consider $V=0$ and equation (\ref{3.8}) reduces to
\begin{equation}\label{eq11}
\frac{D^{4}x}{dt^{4}}+R\left(\frac{D^2x}{dt^2},\frac{dx}{dt}\right)\frac{dx}{dt}- \sigma
\frac{D^2x}{dt^2} = 0
\end{equation} which gives the so called \textit{cubic polynomials in tension} on Riemannian manifolds \cite{SCC00}, that is, smooth trajectories on $M$, given by the extremals among $\Omega$ of the action functional \begin{equation}\label{cubictension}
J(x)=\int_{0}^{T}\frac{1}{2}\left(\Big{\|}\frac{D^2x}{dt^2}(t)\Big{\|}^2+
 \sigma \Big{\|}\frac{ dx}{dt}(t)\Big{\|}^2\right)dt.
\end{equation}
When the parameter $\sigma$ is zero, these curves are \textit{Riemannian cubic polynomials} \cite{CroSil:91}, \cite{Noa:89} and, for nonzero values of $\sigma$, as  $\sigma$ increase the curves approximate more
precisely the geodesics joining the same points. These curves  have many applications in physics and engineering (see for instance \cite{BlochHussein}, \cite{GarayNoakes2014}).
\end{remark}

\section{Dynamic interpolation for obstacle avoidance problems on Riemannian manifolds}
Now we consider the following problem of dynamic interpolation for obstacle avoidance. We start by studying the general case on a Riemannian manifold and then the case of Lie groups. We illustrate the results with several examples.

\subsection{Dynamic interpolation for obstacle avoidance: The general case}
Consider a set of distinct points $x_0,x_1,\ldots,x_N\in M$ such that each $x_j$ does not intersect $S$, $j=0,\ldots,N$,  and a set of  fixed times $0=T_0<T_1<\ldots<T_{N-1}<T_N=T$. We define the admissible set $\overline{\Omega}$ of $\mathcal{C}^{1}$ curves on $[0,T]$, which are smooth on $[T_{i-1},T_i]$, $i=1,\ldots,N$, and verify the interpolation conditions
\begin{equation}\label{interpolc}
x(T_i)=x_i,\quad  i=1,\ldots,N-1,
\end{equation}
and the boundary conditions
\begin{equation}\label{boundaryc}
x(0)=x_0,\quad x(T)=x_N,\quad \frac{dx}{dt}(0)=v_0 \hbox{ and }\frac{dx}{dt}(T)=v_N.
\end{equation}

The tangent space $T_x\overline{\Omega}$ to the curve $x\in\overline{\Omega}$ is defined to be the vector space of $\mathcal{C}^{1}$ vector fields $X$ on $[0,T]$, which are smooth on $[T_{i-1},T_i]$ and satisfy the  conditions
$X(T_i)=0$, $i=0,\ldots,N$ and
\begin{equation}\label{3.}
\frac{DX}{dt}(0)=0,\quad   \frac{DX}{dt}(T)=0.
\end{equation}

\textbf{Problem 2}: The problem of \textit{dynamic interpolation for obstacle avoidance} consists of minimizing the functional $J$ on $\overline{\Omega}$.

\begin{theorem} \label{t3.2}
A necessary condition for $x$  to be an extremizer of the functional \eqref{3.2} over the class $\overline{\Omega}$ is that $x$ is $\mathcal{C}^{2}$  and verifies the following equation
 \begin{equation} \label{3.8b}
\frac{D^{4}x}{dt^{4}}+R\left(\frac{D^2x}{dt^2},\frac{dx}{dt}\right)\frac{dx}{dt}- \sigma
\frac{D^2x}{dt^2}+ \frac{1}{2}\mbox{grad }V(x) = 0
\end{equation}
on each interval $[T_{i-1},T_i]$, $i=1,\ldots,N$.
\end{theorem}

\textbf{Proof}: Let $\alpha$ be an admissible variation of $x$ with variational vector field
$X \in T_x\overline{\Omega}$. Then
\begin{equation*}
\frac{d}{dr}J(\alpha _{r})=\int_{0}^{T}\left(\Big{\langle}\frac{D}{\partial
r}\frac{D^2 \alpha }{\partial t^2},\frac{D^2 \alpha }{\partial t^2}\Big{\rangle}+\sigma\Big{\langle}\frac{D^2\alpha}{\partial r \partial t},\frac{\partial \alpha }{\partial
t}\Big{\rangle}+ \frac{1}{2}\frac{\partial}{\partial r}V(\alpha)\right)dt.
\end{equation*}

By considering the gradient vector field $\mbox{grad }V$ of the potential function $V$ we have $$ \frac{\partial}{\partial
r}V(\alpha)=\Big{\langle}\frac{\partial \alpha}{\partial r },\mbox{grad }V(\alpha)\Big{\rangle}.$$

By Lemma \ref{lemma} and the previous identity we have
\begin{equation*}
\frac{d}{dr}J(\alpha _{r})=\int_{0}^{T}\left(
\Big{\langle}\frac{D^2}{dt^2}\frac{\partial \alpha}{\partial r},\frac{D^2 \alpha }{\partial
t^2}\Big{\rangle}+\Big{\langle}R\left(\frac{\partial \alpha}{\partial r},\frac{\partial
\alpha}{\partial t}\right)\frac{\partial \alpha}{\partial t}, \frac{D^2 \alpha }{\partial t^2}\Big{\rangle}+\sigma
\Big{\langle}\frac{D^2\alpha}{\partial t \partial r},\frac{\partial \alpha }{\partial
t}\Big{\rangle}+\Big{\langle}\frac{\partial \alpha}{\partial r },\frac{1}{2}\mbox{grad }V(\alpha)\Big{\rangle}\right)dt.
\end{equation*}

Since $\alpha$ is smooth on $[T_{i-1},T_i]$, integrating the first term by parts twice, and the third term once, on each interval,
and applying the property (\ref{curvformula}) of the curvature tensor $R$ to the second term, we obtain
\begin{align*}
\frac{d}{dr}J(\alpha _{r})= &
\sum_{i=1}^{N}\left[\Big{\langle}\frac{D}{\partial t}\frac{\partial \alpha}{\partial r},\frac{D^2 \alpha
}{\partial t^2}\Big{\rangle}-\Big{\langle}\frac{\partial \alpha}{\partial r}, \frac{D^3 \alpha }{\partial t^3}\Big{\rangle}+\sigma\Big{\langle}\frac{\partial \alpha }{\partial r},\frac{\partial \alpha }{\partial
t}\Big{\rangle}\right]_{T_{i-1}^+}^{T_i^-}\\&+\int_{0}^{T}
\Big{\langle}\frac{\partial \alpha}{\partial r},\frac{D^4 \alpha }{\partial
t^4}+R\left(\frac{D^2 \alpha}{\partial t^2},\frac{\partial
\alpha}{\partial t}\right)\frac{\partial \alpha}{\partial t}-\sigma\frac{D^2 \alpha }{\partial t^2}+\frac{1}{2}\mbox{grad }V(\alpha)\Big{\rangle}\, dt.
\end{align*}

Next, by taking $r=0$ in the last equality, we obtain
\begin{align*}
\frac{d}{dr}J(\alpha _{r})\Big{|}_{r=0}=&\sum_{i=1}^{N} \left[\Big{\langle}\frac{DX}{dt},\frac{D^2x}{dt^2}\Big{\rangle}- \Big{\langle}X,\frac{D^{3}x}{dt^{3}}\Big{\rangle}+\sigma\Big{\langle} X,\frac{dx}{dt}\Big{\rangle}\right]_{T_{i-1}^{+}}^{T_i^{-}}\\
&+\int_{0}^{T}\Big{\langle}X,\frac{D^{4}x}{dt^{4}}+R\left(\frac{D^2x}{dt^2},\frac{dx}{dt}\right)\frac{dx}{dt}-\sigma \frac{D^2x}{dt^2}+\frac{1}{2}\mbox{grad } V(x(t))\Big{\rangle}\, dt .
\end{align*}

Since the vector field $X$ is ${\mathcal C}^1$, piecewise smooth
on $[0,T]$, verifies the boundary conditions (\ref{3.}) and the curve $x$
is ${\mathcal C}^1$ on $[0,T]$, it follows that, if $\alpha$ is an admissible variation of $x$ with variational vector field
$X \in T_x\overline{\Omega}$, then \begin{align*}\label{3.7}
\frac{d}{dr}J(\alpha _{r})\Big{|}_{r=0}=&\int_{0}^{T}\Big{\langle}X,\frac{D^{4}x}{dt^{4}}+R\left(\frac{D^2x}{dt^2},\frac{dx}{dt}\right)\frac{dx}{dt}-\sigma \frac{D^2x}{dt^2}+\frac{1}{2}\mbox{grad } V(x)\Big{\rangle}\, dt \\
&+\sum_{i=1}^{N-1} \Big{\langle}\frac{DX(T_{i})}{dt},\frac{D^2x}{dt^2}(T^{+}_{i})-
\frac{D^2x}{dt^2}(T^{-}_{i})\Big{\rangle}.
\end{align*}

Now, assume $x$ is an extremizer of $J$
over $\overline{\Omega}$. Then  $\displaystyle{\frac{d}{dr}J(\alpha_r)\mid_{r=0}=0}$, for each admissible variation $\alpha$  of $x$ with variational vector field
$X \in T_x\overline{\Omega}$.

Let us consider  $X \in T_x\overline{\Omega}$ defined by
$$f\left[\frac{D^{4}x}{dt^{4}}+R\left(\frac{D^2x}{dt^2},\frac{dx}{dt}\right)\frac{dx}{dt}- \sigma
\frac{D^2x}{dt^2}+\frac{1}{2}\mbox{grad }V(x)\right],
$$ where $f$ is a smooth real-valued function on $[0,T]$ verifying
$f(T_i)=f'(T_i)=0$ and $f(t)>0$, for all $t\in(T_{i-1},T_i)$.
So, we have
\begin{equation*}
\displaystyle \frac{d}{dr}J(\alpha_r)\Big{|}_{r=0} =\int_0^Tf\Big{|}\Big{|}\frac{D^{4}x}{dt^{4}}+R\left(\frac{D^2x}{dt^2},\frac{dx}{dt}\right)\frac{dx}{dt}- \sigma
\frac{D^2x}{dt^2}+\frac{1}{2}\mbox{grad }V(x) \Big{|}\Big{|}^2dt
\end{equation*} and since $f(t)>0$ for $t\in (T_{i-1},T_i)$, it follows that
$$\Big{|}\Big{|} \frac{D^{4}x}{dt^{4}}+R\left(\frac{D^2x}{dt^2},\frac{dx}{dt}\right)\frac{dx}{dt}- \sigma
\frac{D^2x}{dt^2}+\frac{1}{2}\mbox{grad }V(x)\Big{|}\Big{|}= 0$$ on $[T_{i-1},T_{i}]$ which leads to the equation (\ref{3.8b}) on each subinterval $[T_{i-1},T_i]$.

Finally, let us choose the vector field $X \in T_x\overline{\Omega}$ so that

$$\frac{DX(T_{i})}{dt}=\frac{D^2x}{dt^2}(T^{-}_{i})-
\frac{D^2x}{dt^2}(T^{+}_{i}),$$

\noindent for $i=1,\ldots, N-1$.
Thus,
\begin{equation*}\frac{d}{dr}J(\alpha_r)\Big{|}_{r=0}=\sum_{i=1}^{N-1}\left(\Big{|}\Big{|}\frac{D^2x}{dt^2}(T^{+}_{i})-
\frac{D^2x}{dt^2}(T^{-}_{i})\Big{|}\Big{|}^2\right)=0
\end{equation*}

\noindent which implies that
$$\displaystyle  \frac{D^2x}{dt^2}(T^{+}_{i})=
\frac{D^2x}{dt^2}(T^{-}_{i})
.$$

\noindent Hence, $x$ is $\mathcal{C}^{2}$
on $[0,T]$.
\quad$\Box$
\begin{remark}
In the absence of obstacles, we consider $V=0$ and equation (\ref{eq11}), defined on each subinterval $[T_{i-1},T_i]$ of $[0,T]$, $i=1,\cdots, N$,   gives rise to the generalization of cubic splines in tension to Riemannian manifolds \cite{SCC00}.
\end{remark}

The next result gives an extension of Theorem \ref{t3.2} for multiple obstacles. The proof does not differ of the one given in Theorem \ref{t3.2} except for the term
concerning the potential function.

Assuming that in the workspace we have $s$ obstacles, the functional \eqref{3.2} becomes in
\begin{equation}\label{severalobstacles}
J(x)=\int_{0}^{T}\frac{1}{2}\left(\Big{\|}\frac{D^2x}{dt^2}(t)\Big{\|}^2+
 \sigma \Big{\|}\frac{ dx}{dt}(t)\Big{\|}^2+\sum_{r=1}^{s}V_{r}(x(t))\right)dt,
\end{equation} where each obstacle is represented by $S_r$ and the artificial potential function $V_{r}$ corresponding to  the obstacle $S_r$ is defined as before, $r=1, \ldots, s$.

\begin{corollary}\label{}

A necessary condition for $x$  to be an extremizer of the functional \eqref{severalobstacles} over the class $\overline{\Omega}$ is that $x$ is $\mathcal{C}^{2}$ and verifies the following equation
 \begin{equation} \label{3.8c}
\frac{D^{4}x}{dt^{4}}+R\left(\frac{D^2x}{dt^2},\frac{dx}{dt}\right)\frac{dx}{dt}- \sigma
\frac{D^2x}{dt^2}+ \frac{1}{2}\sum_{r=1}^{s}\mbox{grad }V_{r}(x) = 0
\end{equation}
on each interval $[T_{i-1},T_i]$, $i=1,\ldots,N$.

\end{corollary}

\begin{remarkth}We would like to point out that it is not guaranteed that the action functional $J$ can achieve a minimum value at an interpolating curve. Indeed, in the work \cite{Heeren} authors find conditions on the Riemannian manifolds  for which cubic splines do  not exist (Lemma 2.15 in \cite{Heeren}), that is, non-existence conditions for the critical paths of the dynamic interpolation problem when the artificial potential is zero everywhere and the elastic parameter $\tau
$ is zero.

In \cite{Giambo} authors study the existence of global minimizers for the variational problem (Problem 1) in complete Riemannian manifolds when the artificial potential is zero everywhere and the elastic parameter $\tau
$ is zero (the critical paths correspond  to Riemannian cubic polynomials). Such a result establishes existence conditions for global minimizers by an understanding of the variational problem as one in a Hilbert manifold setting and using techniques of calculus of variations and global analysis on manifolds.

\end{remarkth}

\subsection{Dynamic interpolation for obstacle avoidance problems on a Lie group}

Now we consider a Lie group $G$ endowed with a left-invariant  Riemannian metric $< \cdot , \cdot >$, with $\mathbb{I}:\mathfrak{g}\times \mathfrak{g}\to\R$ the corresponding inner product on the Lie algebra $\mathfrak{g}$, a positive-definite symmetric bilinear form in $\mathfrak{g}$.  The inner product $\mathbb{I}$ defines  the metric $< \cdot , \cdot >$ completely via left translation (see for instance \cite{bookBullo} pp. 273).

The Levi-Civita connection $\nabla$ induced by $< \cdot , \cdot >$ is an affine left-invariant connection and it is completely determined by its restriction to $\mathfrak{g}$ via left-translations. This restriction, denoted by $\stackrel{\mathfrak{g}}{\nabla}:\mathfrak{g}\times\mathfrak{g}\to\mathfrak{g}$,  is given by (see \cite{bookBullo} p. 271)

\begin{equation}\label{restrictedconnection}\stackrel{\mathfrak{g}}{\nabla}_wu= \frac 12 [w,u]-\frac 12 \mathbb{I}^{\sharp}\left(\hbox{ad}_w^* \mathbb{I}^{\flat}(u)+\hbox{ad}_u^* \mathbb{I}^{\flat}(w)\right),\end{equation} where \hbox{ad}$^{*}:\mathfrak{g}\times\mathfrak{g}^{*}\to\mathfrak{g}^{*}$ is the co-adjoint representation of $\mathfrak{g}$ on $\mathfrak{g}^{*}$ and where $\mathbb{I}^{\sharp}:\mathfrak{g}^{*}\to\mathfrak{g}$, $\mathbb{I}^{\flat}:\mathfrak{g}\to\mathfrak{g}^{*}$ are the associated isomorphisms with the inner product $\mathbb{I}$ (see \cite{Boothby} for instance).

We denote by $u_{L}$ the  left-invariant vector field associated  with $u\in\mathfrak{g}$. For the left-invariant vector fields  $u_L$ and $w_L$ , the covariant derivative of  $u_L$ with respect to $w_L$  is given by $\nabla_{w_L}u_L=(\stackrel{\mathfrak{g}}{\nabla}_wu)_L$, for each $u,v\in\mathfrak{g}$.

Let $x:I\subset\mathbb{R}\to G$ be a smooth curve on $G$. The body velocity of $x$ is the curve $v:I\subset\mathbb{R}\to\mathfrak{g}$ defined by $\displaystyle{v(t)=T_{x(t)}L_{x(t)^{-1}}\left(\frac{dx}{dt}(t)\right)}$.

Let $\{e_1,\ldots,e_n\}$ be a basis of  $\mathfrak{g}$. The body velocity of $x$ on the given basis is described  by $\displaystyle{v=\sum_{i=1}^n v_i e_i}$, where $v_1,\ldots, v_n$ are  the so-called pseudo-velocities
of the curve $x$ with respect to the given basis. The velocity vector can be written in terms of the pseudo-velocities as follows.
\begin{equation}\label{admissibility}\frac{dx}{dt}(t)=T_eL_{x(t)}v(t)=\sum_{i=1}^n v_i(t) (e_i)_L(x(t)).\end{equation}

When the body velocity is interpreted as a control on the Lie algebra, equations \eqref{admissibility} give rise to the so called left-invariant control systems discussed in \cite{Noami2}. Therefore our analysis also includes this class of kinematic control systems.

To write the equations determining necessary conditions for optimality, we must use the following formulas (see \cite{Altafini}, Section $7$ for more details).
\begin{align}
&\frac{D^{2}x}{dt^{2}}=T_eL_{x}\Big(v^{\prime}+\stackrel{\mathfrak{g}}{\nabla}_vv\Big),\label{eqq2}\\
&\frac{D^{3}x}{dt^{3}}=T_eL_{x}\Big(v^{\prime \prime}+\stackrel{\mathfrak{g}}{\nabla}_{v^{\prime}}v+2 \stackrel{\mathfrak{g}}{\nabla}_vv^{\prime}+\stackrel{\mathfrak{g}}{\nabla}_v\stackrel{\mathfrak{g}}{\nabla}_vv\Big),\label{eqq3}\\
&\frac{D^{4}x}{dt^{4}}=T_eL_{x}\left(v'''+\stackrel{\mathfrak{g}}{\nabla}_{v''}v+3\stackrel{\mathfrak{g}}{\nabla}_{v^{\prime}}v^{\prime}+3 \stackrel{\mathfrak{g}}{\nabla}_vv''+\stackrel{\mathfrak{g}}{\nabla}_{v^{\prime}}\stackrel{\mathfrak{g}}{\nabla}_vv+2 \stackrel{\mathfrak{g}}{\nabla}_v\stackrel{\mathfrak{g}}{\nabla}_{v^{\prime}}v+3\stackrel{\mathfrak{g}}{\nabla}_v^2v^{\prime}+\stackrel{\mathfrak{g}}{\nabla}^{3}_vv\right),\label{eqq4}\\
&R\left(\frac{D^{2}x}{dt^{2}},\frac{dx}{dt}\right)\frac{dx}{dt}=T_eL_{x}\left(\mathfrak{R}(v^{\prime},v)v+\mathfrak{R}(\stackrel{\mathfrak{g}}{\nabla}_vv,v)v\right)\label{eqq5},\end{align} where $\mathfrak{R}$ denotes the restriction of the curvature tensor to $\mathfrak{g}$.

Using equations \eqref{eqq2}-\eqref{eqq5} and Theorem \ref{t3.2} we obtain the following result.
\begin{corollary}\label{corollary2}
 The equations giving rise to first order necessary conditions for optimality in the problem $2$ defined on a Lie group $G$ are
 \begin{align*}
0=&v'''+\stackrel{\mathfrak{g}}{\nabla}_{v''}v+3\stackrel{\mathfrak{g}}{\nabla}_{v^{\prime}}v^{\prime}+3 \stackrel{\mathfrak{g}}{\nabla}_vv''
+\stackrel{\mathfrak{g}}{\nabla}_{v^{\prime}}\stackrel{\mathfrak{g}}{\nabla}_vv+2\stackrel{\mathfrak{g}}{\nabla}_v\stackrel{\mathfrak{g}}{\nabla}_{v^{\prime}}v+
3\stackrel{\mathfrak{g}}{\nabla}_v^2v^{\prime}+\stackrel{\mathfrak{g}}{\nabla}_v^3v+\mathfrak{R}(v^{\prime},v)v\\&-\sigma \stackrel{\mathfrak{g}}{\nabla}_vv+\mathfrak{R}(\stackrel{\mathfrak{g}}{\nabla}_vv,v)v
-\sigma v^{\prime}+\frac{1}{2}T_xL_{x^{-1}}(\mbox{grad }V(x))\end{align*}
together with equation \eqref{admissibility}, subject to the interpolation conditions $x(T_i)=x_i$, $i=1,\ldots,N-1$, and boundary conditions $x(0)=x_0,$ $x(T)=x_N$,  $v(0)=T_{x_0}L_{x_0^{-1}}(v_0)$, $v(T)=T_{x_N}L_{x_N^{-1}}(v_N)$.
\end{corollary}
As in the previous subsection, in the presence of $s$ obstacles, the previous equation reads

\begin{align*}
0=&v'''+\stackrel{\mathfrak{g}}{\nabla}_{v''}v+3\stackrel{\mathfrak{g}}{\nabla}_{v^{\prime}}v^{\prime}+3 \stackrel{\mathfrak{g}}{\nabla}_vv''
+\stackrel{\mathfrak{g}}{\nabla}_{v^{\prime}}\stackrel{\mathfrak{g}}{\nabla}_vv+2\stackrel{\mathfrak{g}}{\nabla}_v\stackrel{\mathfrak{g}}{\nabla}_{v^{\prime}}v+
3\stackrel{\mathfrak{g}}{\nabla}_v^2v^{\prime}+\stackrel{\mathfrak{g}}{\nabla}_v^3v+\mathfrak{R}(v^{\prime},v)v\\&-\sigma \stackrel{\mathfrak{g}}{\nabla}_vv+\mathfrak{R}(\stackrel{\mathfrak{g}}{\nabla}_vv,v)v
-\sigma v^{\prime}+\frac{1}{2}\sum_{r=1}^{s}T_xL_{x^{-1}}(\mbox{grad }V_{r}(x)).\end{align*}

\begin{example}{Dynamic interpolation for obstacle avoidance on $SE(3)$.}\label{SE(3)}

Next, as an application of Proposition \ref{corollary2}, we study the variational interpolation problem for the motion of the attitude and translation of a rigid body where the configuration space is the  special Euclidean group $SE(3)$ and a spherical obstacle in the workspace must be avoided. Working on $SE(3)$ we represent the orientation and position of the rigid body in a coordinate free framework . This example corresponds to the dynamic interpolation for obstacle avoidance associated with the dynamics of an aerospace or underwater vehicle  (see for instance \cite{bookBullo} p. 281).

We describe the movement of the rigid body by a curve in $SE(3)$. The  special Euclidean group $SE(3)$ consists of all rigid displacements  in $\mathbb{R}^3$, described by a translation after a rotation. Its elements are the transformations of $\mathbb{R}^3$ of the form $z \mapsto Rz+r$, where $r\in \mathbb{R}^3$ and $R\in SO(3)$.

This group has the structure of the semidirect  product Lie group of $SO(3)$ and $\mathbb{R}^3$. Each rigid displacement can be represented by the element $g=(R, r)$ or, in matrix form, by
$\displaystyle{g=\left(
     \begin{array}{cc}
      R & r\\
       0 & 1 \\
     \end{array}
   \right).}$

The composition law is given by $
(R,r)\cdot(S,s)=(RS,Rs+r)$ with identity element $(I,0)$ and inverse
$g^{-1}=(R^{-1},-R^{-1}r)$.
Note that the composition law corresponds to the usual matrix multiplication if we consider the matrix representation.

The Lie algebra $\mathfrak{se}(3)$ of $SE(3)$ is described by the matrices of the form $\xi=\left(
     \begin{array}{cc}
      A & b\\
       0 & 0 \\
     \end{array}
   \right)$, called twists,
with $A\in \mathfrak{so}(3)$ and $b\in \mathbb{R}^3$. A matrix $A\in\mathfrak{so}(3)$, that is, a skew-symmetric matrix of the form $A=\left(
               \begin{array}{ccc}
                 0 &-a_3  & a_2 \\
                 a_3 & 0 & -a_1 \\
                 -a_2 & a_1 & 0 \\
               \end{array}
             \right)$ can be denoted by  $\widehat{a}$, where $a=(a_1,a_2,a_3)\in \mathbb{R}^3$. We identify the Lie algebra $\mathfrak{se}(3)$ with $\mathbb{R}^6$ via the isomorphism
   $\displaystyle{\left(
     \begin{array}{cc}
     \widehat{a}& b\\
       0 & 0 \\
     \end{array}
   \right)\mapsto (a,b)}$.
   The Lie bracket in $\mathbb{R}^{6}$ is given by
$[(a,b),(c,d)]=(a\times c,a\times d-c\times b).$
The elements $(a,b)$  are called twist coordinates.

The adjoint action is given by
$$\hbox{Ad}_{(R,r)}(a,b)=(Ra,Rb-Ra\times r).$$

We consider the basis $\{e_i\}_{i=1}^{6}$ of $\mathfrak{se}(3)$, represented by the canonical basis  of $\mathbb{R}^6$,
  given by
\begin{equation*}
e_1=\left[{\begin{array}{cccc}
   0 & 0 & 0 & 0 \\
   0 & 0 & -1 & 0 \\
   0 & 1 & 0 & 0 \\
   0 & 0 & 0 & 0 \\
  \end{array} }\right], \quad e_2=\left[{\begin{array}{cccc}
   0 & 0 & 1 & 0 \\
   0 & 0 & 0 & 0 \\
   -1 & 0 & 0 & 0 \\
   0 & 0 & 0 & 0 \\
  \end{array} }\right], \quad  e_3=\left[{\begin{array}{cccc}
   0 & -1 & 0 & 0 \\
   1 & 0 & 0 & 0 \\
   0 & 0 & 0 & 0 \\
   0 & 0 & 0 & 0 \\
  \end{array} }\right],\quad e_4=\left[{\begin{array}{cccc}
   0 & 0 & 0 & 1 \\
   0 & 0 & 0 & 0 \\
   0 & 0 & 0 & 0 \\
   0 & 0 & 0 & 0 \\
  \end{array} }\right],\end{equation*}
  \begin{equation*} e_5=\left[{\begin{array}{cccc}
   0 & 0 & 0 & 0 \\
   0 & 0 & 0 & 1 \\
   0 & 0 & 0 & 0 \\
   0 & 0 & 0 & 0 \\
  \end{array} }\right], \quad e_6=\left[{\begin{array}{cccc}
   0 & 0 & 0 & 0 \\
   0 & 0 & 0 & 0 \\
   0 & 0 & 0 & 1 \\
   0 & 0 & 0 & 0 \\
  \end{array} }\right]
  \end{equation*}
  and  endow $SE(3)$ with the left-invariant metric  defined by the inner product $$ \mathbb{I}=\sum_{i=1}^3J_ie^i\otimes e^i+ m_ie^{i+3}\otimes e^{i+3}$$ where $J_i$ and $m_i$ are the diagonal elements of the matrix defining the dynamics of the rigid body, the inertia moments and masses, with $i=1,2,3$ and $\{e^{i}\}_{i=1}^{6}$, the dual basis of $\{e_i\}_{i=1}^{6}$.

The Levi-Civita connection $\nabla$ induced by $< \cdot , \cdot >$ is completely determined by its restriction to $\mathfrak{se}(3)$ and is given by (see for instance \cite{bookBullo} p. 282)
\begin{equation*}
\stackrel{\mathfrak{se}(3)}{\nabla}_{v} z=\frac 12 \Big(a\times c+J^{-1}(a\times J c+b\times M d+c\times J a+d\times M b), b\times c+a\times d+M^{-1}(a\times M d+c\times M b)\Big)\end{equation*} where $J$ and $M$ are the blocks of the diagonal matrix describing the dynamics of the body, representing the moments of inertia and masses respectively, and  $v=(a,b)$, $z=(c,d) \in \mathbb{R}^6$ are the twist coordinates.

 For simplicity in the exposition, we consider the case when $J=M=I$. Then the formula above for the Levi-Civita connection reduces to
$$\stackrel{\mathfrak{se}(3)}{\nabla}_{v} z=\Big(\frac 12 a\times c,b\times c+a\times d\Big),$$
Using \eqref{eq:CurvatureTensorDefinition} we obtain the  restriction of the curvature tensor  to  $\mathfrak{se}(3)$ given by
$$\mathfrak{R}(v,z)w=\Big(-\frac 14 (a\times c)\times f,0\Big).$$
where as before $v=(a,b)$, $z=(c,d)$ and $w=(h,f) \in \mathbb{R}^6$ represent the twist coordinates.

The motion of the rigid body in space is described by a curve $x$ in $SE(3)$. The body velocity is given by the curve $v$ in $\mathfrak{se}(3)$, described in the basis  $\{e_i\}_{i=1}^{6}$ as
$$v=\sum_{i=1}^3 a_i e_i+\sum_{i=1}^3 b_i e_{i+3}.$$

The first three terms correspond to infinitesimal rotations about the three axes (roll, pitch and yaw) and the later three terms to infinitesimal translations about the three axes.

We consider the potential functions \begin{equation}\label{navfunct2}V_1(R,r)=\frac{\tau_1}{\parallel r\parallel^2-1},\quad V_2(R,r)=\frac{\tau_2}{\parallel r-p\parallel^2-2},\end{equation} $R \in SO(3), r\in \R^3\backslash S^2$, designed
for avoidance of two obstacles with spherical shape,  the first with unit radius centered
at the origin and the second
with radius $\sqrt{2}$ centered at $p=(2,2,2)$. Here $\tau_1,\tau_2\in\mathbb{R}^{+}$ and $||\cdot||$ is the Euclidean norm.

We can rewrite
the potential function $V_1$ as follows

\begin{align}
V_1(R,r) = \frac{\tau_1}{\|\hbox{Ad}_{g^{-1}}e_{1}\|^{2}_{\mathfrak{se}(3)}+\|\hbox{Ad}_{g^{-1}}e_{2}\|^{2}_{\mathfrak{se}(3)}+\|\hbox{Ad}_{g^{-1}}e_{3}\|^{2}_{\mathfrak{se}(3)}-1}, \nonumber
\end{align} with $g\in SO(3)\times (\R^3\backslash S^2)$, where $||\cdot||_{\mathfrak{se}(3)}$ is the norm on $\mathfrak{se}(3)$ defined by the inner product on $\mathfrak{se}(3)$  given by $\langle\langle\xi,\xi\rangle\rangle=\hbox{tr}(\xi^{T}\xi)$, for any $\xi\in\mathfrak{se}(3)$. Hence, the norm $\|\xi\|_{\mathfrak{se}(3)}$ is given by $\|\xi\|_{\mathfrak{se}(3)}=\langle\langle\xi,\xi\rangle\rangle^{1/2} = \sqrt{\hbox{tr}(\xi^{T}\xi)}$, for any $\xi\in\mathfrak{se}(3)$. Similarly we can rewrite the potential function $V_2$ on the Lie group using the Adjoint action.

 A form  of Euler-Poincar\'e equations can be obtained as in \cite{Gu} using the $SO(3)$-invariance of $V_1$ and $V_2$. We will study that approach in future work. Here we  study the dynamics using the representation given by the Lie algebra isomorphism $\mathfrak{se}(3)\simeq\mathbb{R}^6$. The gradient of $V_1$ satisfies
 \begin{align*}T_gL_{g^{-1}}(\mbox{grad } V_1(g))&=T_gL_{g^{-1}}\left(\frac{-\tau_1(2x\partial_x +2y\partial_y+2z\partial_z)}{(\parallel r\parallel^2-1)^2}\right)\\
&=-\frac {2\tau_1}{(\parallel r\parallel^2-1)^2}\left(
                    \begin{array}{cc}
                      R^T & -R^Tr\\
                      0 & 1\\
                    \end{array}
                  \right)\left(
                           \begin{array}{cc}
                             0 & r \\
                             0 & 0\\
                           \end{array}
                         \right)\\
                         &=
 -\frac {2\tau_1}{(\parallel r\parallel^2-1)^2}\left(0,R^Tr\right)\end{align*}
 where here, we have  $x\partial_x +y\partial_y+z\partial_z=\left(
                               \begin{array}{cccc}
                                    0 & 0& 0& x \\
                                   0& 0 & 0& y\\
                                   0& 0 & 0& z\\
                                   0 & 0& 0 &0\\
                                  \end{array}
                                \right)$ and $T_gL_{g^{-1}}$ is the product by $\left(
                    \begin{array}{cc}
                      R^T & -R^Tr\\
                     0 & 1\\
                    \end{array}
                  \right)$ on the left. Similarly, the gradient of $V_2$ satisfies $$T_gL_{g^{-1}}(\mbox{grad } V_2(g))=-\frac {2\tau_2}{(\parallel r-p\parallel^2-2)^2}\left(0,R^Tr\right).$$

By Proposition \ref{corollary2} the necessary conditions for the extremizer in problem $2$ are determined by the equations
\begin{align*}\label{eqbody}
a'''&= a''\times a + \sigma a',\\
b'''&=3b''\times a + \sigma b'+3b'\times a'-3(b'\times a)\times a+ b\times a''+\sigma b\times a-3(b\times a')\times a+b\times (a'\times a) -(b\times a)\times a \\&+\frac {\tau_1}{(\parallel r\parallel^2-1)^2}R^Tr+\frac {\tau_2}{(\parallel r-(2,2,2)\parallel^2-2)^2}R^Tr,\end{align*}
together with the equations
\begin{equation}\label{bodycoord}
R'=R\widehat{a}, \; r'=Rb,
\end{equation} the interpolation conditions $(R(T_i),r(T_i))=(R_i,r_i)$ and the boundary conditions $(R(0),r(0))=(R_0,r_0)$,  $(R(T),r(T))=(R_N,r_N)$, $(a(0), b(0))=(a_0,b_0)$, $(a(T),b(T))=(a_N,b_N)$, where $(a_0,b_0)$ and $(a_N,b_N)$ are the twist coordinates of the body velocity of the curve $(R,r)$ at $0$ and $T$.

Note that, in the absence of obstacles, the extremals  reduce to the  cubic splines in tension on $SE(3)$ \cite{SCC00} given by the following equations.
\begin{align*}
a'''&= a''\times a + \sigma a',\\
b'''&=3b''\times a + \sigma b'+3b'\times a'-3(b'\times a)\times a+ b\times a''+\sigma b\times a-3(b\times a')\times a+b\times (a'\times a) -(b\times a)\times a.\end{align*}

\begin{figure}[h]
\includegraphics[width=4.5cm]{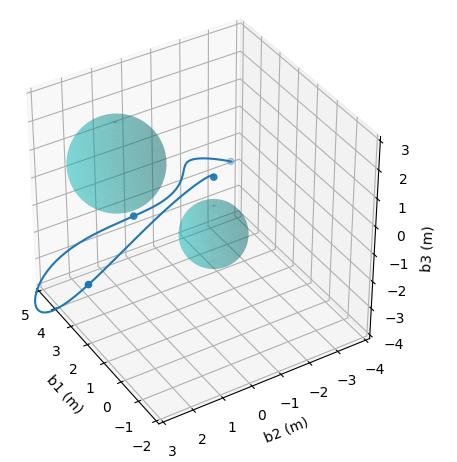}\label{figure1}\includegraphics[width=4.5cm]{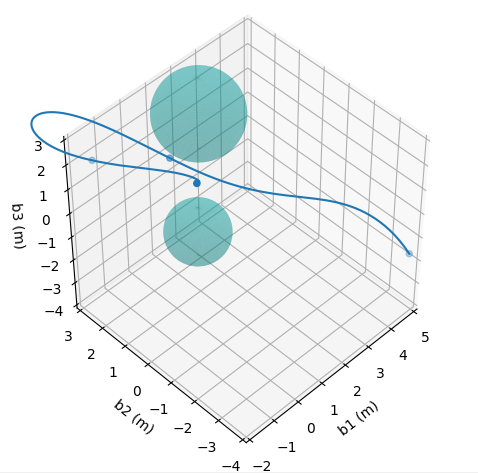}
\caption{Smooth trajectory generation for two obstacle avoidance problem and two interpolation points, given boundary points.}
\end{figure}

In Figure $2$ we show a simulation of our method. A shooting method and a symplectic Euler discretization with $h=0.01$ are used to simulate the boundary value problem. The  curve represents the optimal trajectory interpolating the prescribed points and boundary values. One interpolation point has been taken to be close to one obstacle and between the two prescribed obstacles.

The parameters for the trajectory used are $\sigma= 0.5$, $\tau_1 = 1.7$ and $\tau_2 = 1.1$. Boundary condition are given by:
\[
R_0=
  \begin{bmatrix}
    1 & 0 & 0  \\
    0& 1 & 0 \\
    0 & 0 & 1
  \end{bmatrix}, \quad R_T=
  \begin{bmatrix}
    \frac{\sqrt{2+\sqrt{2}}}{2} & 0 & -\frac{\sqrt{2-\sqrt{2}}}{2}  \\
    0& 1 & 0 \\
  \frac{\sqrt{2-\sqrt{2}}}{2} & 0 & \frac{\sqrt{2+\sqrt{2}}}{2}
  \end{bmatrix},
\] $r_0=(-1,1,0)$, $r_T=\frac{1}{2}(-5,\frac{19}{2}, 3)$, $a_0=(0,0,0)$, $a_T=(-1,-1,4)$, $b_0=(0,0,2)$, $b_T=(5,-4,-2)$, initial time $T=0$ and final time $T=1.87$.

The interpolation points are $\gamma(0.74)=(R_1,r_1)$ and $\gamma(1.43)=(R_2,r_2)$, where
\[
R_1=
   \begin{bmatrix}
    \frac{\sqrt{3}}{2} & 0 & -\frac{1}{2}  \\
    0& 1 & 0 \\
    \frac{1}{2} & 0 & \frac{\sqrt{3}}{2}
  \end{bmatrix}, \quad R_2=
  \begin{bmatrix}
    -\frac{1}{2}& 0 & -\frac{\sqrt{3}}{2}  \\
    0& 1 & 0 \\
    \frac{\sqrt{3}}{2}  & 0 &-\frac{1}{2}
  \end{bmatrix},
\] $r_1=(-1,2,1)$ and $r_2=\frac{1}{2}(-3,11,1)$.
\end{example}

\subsection{Dynamic Interpolation for obstacle avoidance problems on compact and connected Lie groups}
Next, we derive the equations for the dynamic interpolation problem obtained in the previous subsection in the case when the Lie group is compact and connected.

Assume $G$ is a connected and compact Lie group. Therefore $G$ is endowed with a bi-invariant Riemannian metric that makes $G$ a complete Riemannian manifold. In this context the Riemannian distance between two points in $G$ can be defined by means of the Riemannian exponential on $G$, that is,
$$d(g,h)=\|\mbox{exp}_h^{-1}g\|=d(g,h)\quad g,h\in G.$$
We need to guarantee that the exponential map $\mbox{exp}_h$ is a local diffeomorphism, so  we assume that the point $g$ must belong to a convex open  ball around $h$.   If we consider the geodesic from $g$ to $h$ given by
$\gamma_{g,h}(s)= \mbox{exp}_g(s\, \mbox{exp}_g^{-1}h)$, $s\in[0,1]$, then, because $\displaystyle{\Big{\|}\frac{d \gamma_{g,h}}{d s}(s)\Big{\|}}$ is independent of $s$, we can write $$d^2(g,h)=\int_{0}^{1}\Big{\|}\frac{d \gamma_{g,h}}{d s}(s)\Big{\|}^2\, ds.$$
The obstacle is represented by an element $h$ in $G$ and the artificial potential function, used to avoid the obstacle, is defined by
$$V_h(g)=\frac{\tau}{d^2(h,g)}.$$

If we consider a map $\alpha: r\to \alpha(r)$ verifying $\alpha(0)=g$ and  the family of geodesics from $g$ to $\alpha(r)$ given by
$\gamma(s,r)=\mbox{exp}_h(s\, \mbox{exp}_h^{-1}\alpha(r)),$ then we have
$$\displaystyle \frac{d}{d r}d^2(g,\alpha)=\Big{\langle}\frac{\partial \gamma}{\partial r },\frac{\partial \gamma}{\partial
s}\Big{|}_{s=1} \Big{\rangle}=-\Big{\langle}\frac{d \alpha}{dr },\mbox{exp}_{\alpha}^{-1}h\Big{\rangle}$$
and we obtain the expression of the gradient vector field  as follows:
\begin{equation}\label{gradV}\mbox{grad }V_h(g)=\frac{\tau}{d^4(h,g)}\exp_{g}^{-1}h.\end{equation}

The Levi-Civita connection and the curvature tensor of $G$, when restricted to the Lie algebra $\mathfrak{g}$ of $G$, are defined by  \begin{align}
&\stackrel{\mathfrak{g}}{\nabla}_wu= \frac 12 [w,u],\label{eqq2-2}\\
&\mathfrak{R}(w,u)z=-\frac 14[[w,u],z].\label{eqq4-2}\end{align}

Consider as before the body velocity $v$ of a curve $x$ on the Lie group $G$ with respect to a basis $\{e_1,\ldots,e_n\}$ of the Lie algebra $\mathfrak{g}$.
Using \eqref{eqq2-2}-\eqref{eqq4-2}, equation \eqref{gradV} and Corollary \ref{corollary2}, the equations describing first-order optimality conditions for  Problem $2$ on a connected and compact Lie group are obtained as follows:

\begin{corollary}\label{corollary3}
Let $x$ be a $C^2$-curve on a connected and compact Lie group $G$ with body velocity $v$ with respect to the basis $\{e_1,\ldots,e_n\}$ of $\mathfrak{g}$. If the curve $x$  is an extremizer of the functional \eqref{3.2} over the class $\overline{\Omega}$, then $x$ verifies
 \begin{equation} \label{completeG}
v'''-\sigma v'+[v,v'']+\frac{\tau}{2||\exp_{x}^{-1}h||^4}T_xL_{x^{-1}}(\exp_{x}^{-1}h) \equiv 0
\end{equation}
on each interval $[T_{i-1},T_i]$, $i=1,\ldots,N$ where $\exp$ is the exponential map on $G$ and $h\in G$.
\end{corollary}

\begin{example}{Dynamic interpolation for obstacle avoidance problem on $SO(3)$.}

Motivated by the fact that obstacle avoidance problems defined on the special orthogonal group $SO(3)$ are often used to avoid certain pointing directions/orientations (for example avoiding pointing an optical instrument at the Sun) we consider the following interpolation obstacle avoidance problem. Consider a rigid body where the configuration space is the Lie group $G=SO(3)$. The Lie algebra $\mathfrak{so}(3)$ is given by the set of $3\times 3$ skew-symmetric matrices. It is well known that via the hat operator defined in Example \ref{SE(3)} one has the identification $\mathfrak{so}(3)\equiv \R^3$.

For simplicity in the exposition we consider the case of a symmetric rigid body, so $SO(3)$ is endowed with the bi-invariant metric  defined by the Euclidean inner product in $\mathbb{R}^3$. Then the formulas  \eqref{eqq2-2} for the restriction of  the Levi-Civita connection reduces to
$\displaystyle{\stackrel{\mathfrak{so}(3)}{\nabla}_{w} u=\frac{1}{2}w\times u}$
and, using \eqref{eqq4-2}, the  restriction of the curvature tensor  to  $\mathfrak{so}(3)$ is defined by
$\displaystyle{\mathfrak{R}(w,u)z=-\frac 14 (w\times u)\times z}$
where $w,u,z \in \mathbb{R}^3$.

The motion of the rigid body in space is described by a curve $t\to R(t)$ in $SO(3)$. The columns of the
matrix $R(t)$ represent the directions of the principal axis of the
body at time $t$ with respect to some reference system.  The body angular velocity is given by the curve $t\to v(t)$ on $\mathfrak{so}(3)$, described as $v(t)=(v_1(t),v_2(t),v(t))$ in $\mathbb{R}^3$. For the obstacle avoidance problem we consider the navigation functions $V:SO(3)\to\mathbb{R}$ given by \begin{equation*}\label{navfunct3}V_{Q}(R)=\frac{\tau_i}{\|\mbox{exp}_{Q}^{-1}R\|^2}.
\end{equation*} representing a repulsive potential function to avoid the  obstacle  $Q\in SO(3)$, with $\tau\in\mathbb{R}^{+}$, $R\in SO(3)$ and $\exp$ representing the exponential map on $SO(3)$.

Given that $$T_RL_{R^T}(\mbox{grad } V_{Q}(R))=\frac{\tau_i}{\|\mbox{exp}_{Q}^{-1}R\|^4}\exp_I^{-1}({R}^{T}Q),$$ the necessary conditions for optimality are determined by the  equation

$$v'''= v''\times v + \sigma v'+\frac{\tau_1}{2\|\mbox{exp}_{Q}^{-1}R\|^4}\exp_I^{-1}({R}^{T}Q).$$ together with the equation $R'=Rv$, the interpolation points $R(T_i)=R_i$ and the boundary conditions $R(0)=R_0$, $R(T)=R_N$, $v(0)=v_0$ and $v(T)=v_N$.

Note that in the absence of obstacles, the extremals  reduce to the cubic splines in tension on $SO(3)$, \cite{SCC00} where the  equations are given by solutions of the equation $ v'''= v''\times v + \sigma v'$ \end{example}

\section{Sub-Riemannian dynamic interpolation for obstacle avoidance}

Next, we extend our analysis to the sub-Riemannian setting, that is, we consider that the velocity vector field $\displaystyle{\frac{dx}{dt}}$ lies on some
distribution $\mathcal{D}$ on $M$. We assume that $\mathcal{D}$  is a constant dimensional  nonintegrable distribution and that there exist $k$ linearly independent one-forms \linebreak  $\omega _1,\ldots, \omega _k$, such  that the codistribution annihilating $\mathcal{D}$  is spanned by $\omega _1,\ldots, \omega _k$ ($k<n$). The constraints on the velocity vector field are defined by
\begin{equation}\label{eqconstraint}
\omega_j\left(\frac{dx}{dt}\right)=\Big{\langle}Y_j,\frac{dx}{dt}\Big{\rangle}=0,  \quad  j=1, \ldots, k,
\end{equation} where $Y_1, \cdots ,Y_k$ are linearly independent vector fields on $M$.

To deal with the constraints we also need to define the tensors $S_{j}$, $(S_{j})_{x}:T_{x}M\to T_{x}M$, given by
\begin{equation}\label{eqconstraint2}d\omega_j(u,z)=<(S_{j})_{x}(u),z>=-<(S_{j})_{x}(z),u>, u,z \in T_xM.
\end{equation}

\textbf{Problem 3}: The sub-Riemannian dynamic interpolation problem for obstacle avoidance consists of minimizing the functional $J$ defined in \eqref{3.1} on $\overline{\Omega}$ with the additional constraints \eqref{eqconstraint}.

This type of problem, in the absence of obstacles, was studied in Bloch and Crouch \cite{BlCr} and Crouch and Silva Leite \cite{CroSil:95}.

 We derive optimality conditions for this sub-Riemannian problem, by extending our previous analysis for the general case following the result of Bloch and Crouch \cite{BlCr}, \cite{blochcrouch}.
\begin{theorem} \label{T4}
A necessary condition for $x\in \overline{\Omega} $ to give a normal extremum for problem 3 is that $x$ is  ${\mathcal C}^2$ and there exist smooth
functions $\lambda _j$, $j=1,\cdots ,k$ (the Lagrange multipliers) such
that the following equation holds
\begin{equation*}\label{eq:23}
0=\frac{D^{4}x}{dt^{4}}+R\left(\frac{D^2x}{dt^2},\frac{dx}{dt}\right)\frac{dx}{dt}- \sigma
\frac{D^2x}{dt^2}+\frac{1}{2}\mbox{grad }V(x)-\sum_{j=1}^k\lambda _j^{\prime}Y_j-\sum_{j=1}^k\lambda _jS_j\left(\frac{dx}{dt}\right),
\end{equation*}
on each interval $[T_{i-1},T_i]$, $i=1,\ldots,N$, together with $\Big{\langle}Y_j,\frac{dx}{dt}\Big{\rangle}=0, \; j=1,\cdots ,k$.

\end{theorem}

\textbf{Proof}: Consider the extended functional
\begin{align*}
\widetilde{J}(x)=&\frac{1}{2}\int_{0}^{T}\left(\Big{\|}\frac{D^2x}{dt^2}(t)\Big{\|}^2+
 \sigma \Big{\|}\frac{dx}{dt}(t)\Big{\|}^2+V(x(t))+\sum_{j=1}^{k}\lambda_j\Big{\langle}Y_j,\frac{dx}{dt}\Big{\rangle}\right)dt.\end{align*}

We derive necessary conditions for existence of normal extremizers by studying the equation $$\frac{d}{dr}\tilde{J}(\alpha_r)\Big{|}_{r=0}=0$$ for $\alpha$ an admissible variation of $x$ with variational vector field $X \in T_x\Omega$ and $\lambda_j$ the Lagrange multipliers.

Taking into account the proof of Theorem \ref{t3.2} we only need to study the influence of variations in the term $\displaystyle{\sum_{j=1}^{k}\lambda_j\Big{\langle}Y_j,\frac{dx}{dt}\Big{\rangle}}$ where the vector fields $Y_{j}$ on $M$ are determined by $\omega_{j}(Z)=\langle Y_{j}, Z\rangle$, $j=1,\ldots,k$ for each vector field  $Z$ on $M$. Therefore, $\displaystyle{\frac{d}{dr}\tilde{J}(\alpha_r)\Big{|}_{r=0}}$ must have two additional terms compared with $\displaystyle{\frac{d}{dr}J(\alpha_r)\Big{|}_{r=0}}$. Those terms are $$\sum_{j=1}^{k}\lambda_{j}\Big{\langle}\nabla_{\frac{\partial\alpha}{\partial r}}Y_{j},\frac{\partial\alpha}{\partial t}\Big{\rangle}+\sum_{j=1}^{k}\lambda_{j}\Big{\langle}Y_{j},\frac{D^2\alpha}{\partial t \partial r}\Big{\rangle}.$$

After integration by parts in the second term and evaluating at $r=0$, the integrand can be re-written with the additional terms $$\sum_{j=1}^{k}\lambda_j\langle\nabla_{X}Y_{j},\frac{dx}{dt}\rangle-\sum_{j=1}^{k}\lambda_{j}'\langle Y_{j},X\rangle-\lambda_j\sum_{j=1}^{k}\langle\frac{DY_{j}}{dt},X\rangle.$$ Using the identity (\ref{difomega})  the new terms compared with the ones provided by Theorem \ref{t3.2}  which give rise to optimality conditions for $x$ to be a normal extremizer in this sub-Riemannian problem are: $$-\sum_{j=1}^{k}\lambda_{j}d\omega_j\left(\frac{dx}{dt},X\right)-\sum_{j=1}^{k}\lambda_j'\Big{\langle}Y_j,X\Big{\rangle}.$$ Using the fact that $d\omega_j\left(\frac{dx}{dt},X\right)=\Big{\langle}S_j\left(\frac{dx}{dt}\right),X\Big{\rangle}$ the result follows. \quad$\Box$

\begin{remarkth}
The introduction of constraints in the velocities for the collision avoidance variational problem causes difficulties in the study of both geometrical and analytical aspects, as remarked in \cite{Giambo2}, and leads to sophisticated situations as when abnormal minimizers appear. As far as we know there is no definitive result  yet regarding  existence and regularity for minimizers in this sub-Riemannian variational problem. It would be very interesting to explore the geometrical and analytical aspects for the existence of minimizers in the collision avoidance problem under constraints in the velocities.
\end{remarkth}

\begin{corollary}\label{cor4.2}
Any abnormal extremizer for the sub-Riemannian dynamic interpolation for obstacle avoidance satisfy $$\sum_{j=1}^k\lambda _j^{\prime}Y_j+\sum_{j=1}^k\lambda _jS_j\left(\frac{dx}{dt}\right)=0,$$
 where $\lambda_j$, $j=1,\ldots,k$ are not all identically zero.\end{corollary}

The following corollaries are direct consequences of the results presented in Section $3$ for the sub-Riemannian problem by straightforward modifications in the proof of Theorem \ref{T4}.

\begin{corollary}\label{cor4.3} If the number of obstacles on $M$ is $s$, located at the points $q_{r}\in M$, $r=1,\ldots s$,  a necessary condition for $x\in\overline{\Omega}$ to be a normal extremizer  for the sub-Riemannian dynamic interpolation for obstacle avoidance is that $x$ is  ${\mathcal C}^2$ and there exist smooth
functions $\lambda _j$, $j=1,\cdots ,k$ (the Lagrange multipliers) such
that the following equation holds
\begin{equation*}
0=\frac{D^{4}x}{dt^{4}}+R\left(\frac{D^2x}{dt^2},\frac{dx}{dt}\right)\frac{dx}{dt}- \sigma
\frac{D^2x}{dt^2}+\frac{1}{2}\sum_{r=1}^{s}\mbox{grad }V_{r}(x)-\sum_{j=1}^k\lambda _j^{\prime}Y_j-\sum_{j=1}^k\lambda _jS_j\left(\frac{dx}{dt}\right),
\end{equation*} on each interval $[T_{i-1},T_i]$, $i=1,\ldots,N$,  together with $ \Big{\langle}Y_j,\frac{dx}{dt}\Big{\rangle}=0, \; j=1,\cdots ,k$.

\end{corollary}

Now we consider the problem on a Lie group $G$ as we did in section 3.2. We suppose that the constraints on the velocity vector field are defined by a left-invariant distribution $\mathcal{D}$. Let $\{e_1,\ldots,e_n\}$ be an orthogonal  basis of the Lie algebra $\mathfrak{g}$ in such a way that the constraints are given by the left-invariant one-forms $\omega_j$ associated with $e_j$, $j=1,\ldots,k$. We have, as before, \begin{equation}\label{eqconstraintLie}
\omega_j\left(\frac{dx}{dt}\right)=\Big{\langle}Y_j,\frac{dx}{dt}\Big{\rangle}=0, \quad j=1,\ldots,k,
\end{equation} but now $Y_1, \cdots ,Y_k,\cdots , Y_n$ form a basis of orthogonal left-invariant vector fields associated with the elements of the basis of $\mathfrak{g}$. Moreover, since the basis of $\mathfrak{g}$ is orthogonal, the
 vector fields $Y_{k+1},\cdots , Y_n$ span the distribution $\mathcal{D}$. Furthermore, the 2-forms $d\omega_j$ are left-invariant and there exist linear maps
 $\mathfrak{S}_{j}:\mathfrak{g}\to\mathfrak{g}$ such that $d\omega_j(z,w)=\mathbb{I}(\mathfrak{S}_{j}(z),w)$, $z,w\in \mathfrak{g}$. These maps can be expressed by  $\mathfrak{S}_{j}=T_{x}L_{x^{-1}}\circ(S_j)_x\circ T_{e}L_{x}$ with $(S_j)_x$ as in \eqref{eqconstraint2}.

\begin{corollary}\label{cor4.4}

A necessary condition for $x\in\overline{\Omega}$  to be a normal extremizer for the problem $3$ in the Lie group $G$ is that $x$ is ${\mathcal C}^2$ and there exist smooth functions $\lambda _j$, $j=1,\cdots ,k$ (the Lagrange multipliers) such that  the following equation holds
\begin{align*}
0=&v'''+\stackrel{\mathfrak{g}}{\nabla}_{v''}v+3\stackrel{\mathfrak{g}}{\nabla}_{v^{\prime}}v^{\prime}+3 \stackrel{\mathfrak{g}}{\nabla}_vv''
+\stackrel{\mathfrak{g}}{\nabla}_{v^{\prime}}\stackrel{\mathfrak{g}}{\nabla}_vv+2\stackrel{\mathfrak{g}}{\nabla}_v\stackrel{\mathfrak{g}}{\nabla}_{v^{\prime}}v+
3\stackrel{\mathfrak{g}}{\nabla}_v^2v^{\prime}+\stackrel{\mathfrak{g}}{\nabla}_v^3v+\mathfrak{R}(v^{\prime},v)v-\sigma \stackrel{\mathfrak{g}}{\nabla}_vv+\mathfrak{R}(\stackrel{\mathfrak{g}}{\nabla}_vv,v)v
-\sigma v^{\prime}\\&+\frac{1}{2}\sum_{r=1}^{s}T_xL_{x^{-1}}(\mbox{grad }V_{r}(x))
-\sum_{j=1}^k\lambda _j^{\prime}e_j-\sum_{j=1}^k\lambda _j\mathfrak{S}_{j}\left(v\right)
\end{align*}
on each interval $[T_{i-1},T_i]$, $i=1,\ldots,N$, together with  the equation \eqref{admissibility} and the constraints  $ v_j=0, \; j=1,\ldots,k$, subject to boundary conditions $x(0)=x_0, x(T)=x_N, v(0)=T_{x_0}L_{x_0^{-1}}(v_0)$, $v(T)=T_{x_N}L_{x_N^{-1}}(v_N)$, and the interpolation conditions $x(T_i)=x_i$, $i=1,\ldots,N-1$.

\end{corollary}

\begin{example}{Dynamic interpolation for obstacle avoidance  of a unicycle.}

We study the motion planning of a unicycle with one obstacle in the workspace.
The unicycle is a homogeneous disk on a horizontal plane and it is equivalent to  a knife edge on the plane \cite{Bl, bookBullo}.
The configuration of the unicycle at any given time is completely determined by an element of the  special Euclidean group $\mathrm{SE}(2)$.

 The elements of $SE(2)$ can be described by transformations of $\mathbb{R}^2$ of the form $z \mapsto Rz+r$, where $r\in \mathbb{R}^2$ and $R\in SO(2)$. The transformations can be represented by $(R, r)$, where
$$R=\left(
     \begin{array}{cc}
       \cos \theta & -\sin \theta\\
       \sin \theta & \cos \theta \\
     \end{array}
   \right)\quad \mbox{and} \quad r=\left(
     \begin{array}{cc}
       x\\
     y\\
     \end{array}
   \right)$$
or, for the sake of simplicity, by the matrix
$\left(
     \begin{array}{cc}
      R & r\\
       0 & 1 \\
     \end{array}
   \right).$

The composition law is defined by
$(R,r)\cdot(S,s)=(RS,Rs+r)$
with identity element $(I,0)$ and inverse
$(R,r)^{-1}=(R^{-1},-R^{-1}r)$. The  special Euclidean group $SE(2)$ has the structure of the semidirect  product Lie group of $SO(2)$ and $\mathbb{R}^2$.

The Lie algebra $\mathfrak{se}(2)$ of $SE(2)$ is determined by
$$\mathfrak{se}(2)=\Big{\{}\left(
     \begin{array}{cc}
      A & b\\
       0 & 0 \\
     \end{array}
   \right): A\in \mathfrak{so}(2) \hbox{ and } b\in \mathbb{R}^2\Big{\}}.$$
For simplicity, we write $A=-a \mathbb{J} $, $a\in\mathbb{R}$, where $ \mathbb{J} =\left(
     \begin{array}{cc}
      0& 1\\
       -1 &0 \\
     \end{array}
   \right)$ and we identify the Lie algebra $\mathfrak{se}(2)$ with $\mathbb{R}^3$ via the isomorphism
   $\displaystyle{\left(
     \begin{array}{cc}
     -a \mathbb{J} & b\\
       0 & 0 \\
     \end{array}
   \right)\mapsto (a,b)}$.

   The Lie bracket in $\mathbb{R}^{3}$ is given by
$[(a,b),(c,d)]=(0,-a \mathbb{J} d+c \mathbb{J} b).$
The  basis $\{e_i\}_{i=1}^{3}$ of  $\mathfrak{se}(2)$ represented by the canonical basis of $\mathbb{R}^3$ verifies
$[e_1,e_2]=e_3$, $[e_2,e_3]=0$, $[e_3,e_1]=e_2.$
We  endow $SE(2)$ with the left-invariant metric  defined by the inner product $$ \mathbb{I}=Je^1\otimes e^1+ m(e^{2}\otimes e^{2}+e^{3}\otimes e^{3}),$$ where $m$ and  $J$ are the  mass  of the body and its inertia moment about the center of mass and and $\{e^{i}\}_{i=1}^{3}$ is the dual basis of $\{e_i\}_{i=1}^{3}$.

The Levi-Civita connection $\nabla$ induced by $< \cdot , \cdot >$ is defined by its restriction to $\mathfrak{se}(2)$ given by$\stackrel{\mathfrak{se}(2)}{\nabla}:\mathfrak{se}(2)\times\mathfrak{se}(2)\to\mathfrak{se}(2)$ and given by $$\stackrel{\mathfrak{se}(2)}{\nabla}_z w=-a(d_2e_2-d_1e_3)=\left(\begin{array}{c}
                                   0 \\
                                   -ad_2 \\
                                   ad_1
                                 \end{array}\right)=(0,-a \mathbb{J}d)
,$$ where $v=(a,b)=(a,b_1,b_2)$ and $w=(c,d)=(c,d_1,d_2)$ are the representative elements of $\mathfrak{se}(2)$ in $\mathbb{R}^3$  (see \cite{bookBullo} p. 279).
The curvature tensor is zero.

We consider the  potential functions $V_1(R,r)$ and $V_2(R,r)$ given by
\begin{equation*}\label{navfunct}V_1(R,r)=\frac{\tau_1}{\|r\|^2-1},\quad  V_2(R,r)=\frac{\tau_2}{\|r-p\|^2-2}.\end{equation*}
$V_1$ and $V_2$ are introduced to avoid two obstacles with circular shape in the $xy$-plane.  The first has unit radius and is
centered at the origin. The second has  radius $\sqrt{2}$  and is
centered at $p=(2,2)$.  $\tau_1,\tau_2\in\mathbb{R}^{+}$ and $||\cdot||$ is the Euclidean norm.

The knife edge constraint is defined by the  one-form $\omega=\sin \theta dx-\cos \theta dy$ whose associated vector field with respect to the Riemannian metric is
 $$Y= \frac 1 m \left(\sin \theta \frac{\partial}{\partial x}-\cos \theta \frac{\partial}{\partial y}\right).$$
 Note that  the tensor $S$ associated with the knife edge constraint is defined by (\ref{eqconstraint2}) and satisfies
$$S(U)=-\frac 1J(u_2\cos \theta +u_3\sin \theta)\frac{\partial}{\partial \theta}+\frac {u_1}m(\cos \theta  \frac{\partial}{\partial x}+\sin \theta \frac{\partial}{\partial y}),$$
for each  vector field on $SE(2)$ denoted by $\displaystyle U=u_1 \frac{\partial}{\partial \theta}+ u_2 \frac{\partial}{\partial x}+u_3 \frac{\partial}{\partial y}$.  Here, we think of $S(U)(x)$ as $S_x(U(x))$,  $x\in M$, where $S_x:T_{x}M\to T_{x}M$, $x\in M$,  is the tensor $S$  defined in (26) and $U: M\to TM; x\mapsto U(x)$ is a vector field.

We consider a basis of vector fields $\{Y_1,Y_2,Y_3\}$ defined by
 $\displaystyle{Y_1=\frac{1}{J}\frac{\partial}{\partial \theta}}$, $\displaystyle{Y_2=\frac{\cos \theta}{m} \frac{\partial}{\partial x}+\frac{\sin \theta}{m}\frac{\partial}{\partial y}}$ and $Y_3=Y$.
This  is  the basis of left-invariant vector fields associated with  $\displaystyle \frac 1 J e_1$, $\displaystyle \frac 1 m e_2$ and  $\displaystyle \frac 1 m e_3$.
The  distribution $\mathcal{D}$ spanned by $Y_1$ and $Y_2$ and the  one-form $\omega$ are in the conditions of Corollary \ref{cor4.4}. The map $\mathfrak{S}$ corresponding to the tensor $S$ is given by $\displaystyle \mathfrak{S}(u)=-\frac 1 J u_2e_1+\frac 1 m u_1e_2$, for each $\displaystyle u=\sum_{i=1}^3u_ie_i\in \mathfrak{se}(2)$.

By Corollary \ref{cor4.4} the equations determining necessary conditions for normal extremizers in  problem $3$ are

\

\begin{align*}
a'''&=\sigma a^{\prime} - \frac  1 J \lambda b_1\\
b'''&=3 a{a^{\prime}}b+3a^2b^{\prime}-(a^3-a'')\mathbb{J}b+3a^{\prime}\mathbb{J}b^{\prime}+3a\mathbb{J}b''+\sigma (b^{\prime}-a\mathbb{J}b)+ \frac {\tau_1}{m(\parallel r\parallel^2-1)^2}R^Tr\\&+ \frac {\tau_2}{m(\parallel r-p\parallel^2-2)^2} R^Tr +\frac  1 m \lambda \left(\begin{array}{c}
                                  a \\
                                   0
                                 \end{array}\right)+\frac 1 m\lambda'\left(\begin{array}{c}
                                 0\\
                                  1
                                 \end{array}\right)
\end{align*}
together with  $R^{\prime}=-a\mathbb{J}R, r^{\prime}=Rb$ and
the constraint $b_2=0$, where $b=(b_1,b_2)$.

\begin{figure}[h!]
\includegraphics[width=4.4cm]{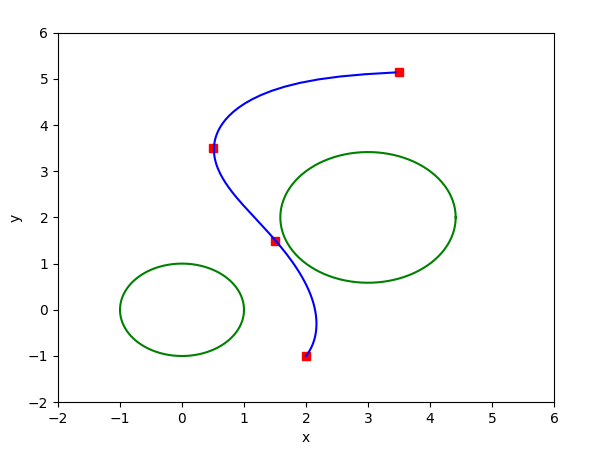}\label{figure1}\includegraphics[width=4.5cm]{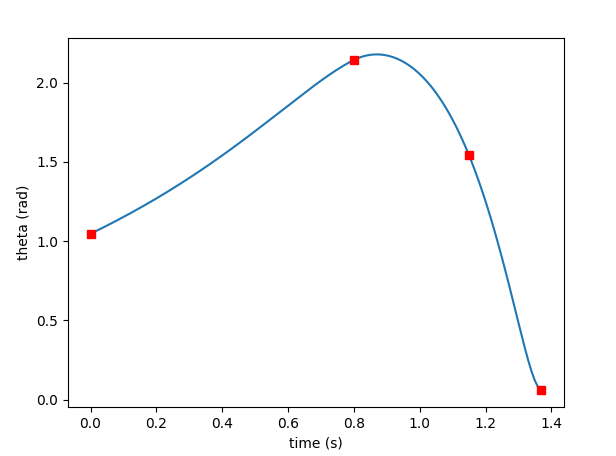}
\caption{Left: Smooth trajectory generation for two obstacle avoidance problem and two interpolation points and given boundary points. Right: evolution of $\theta$ (mod $2\pi$) along the time.}
\end{figure}%

  In Figures \eqref{figure1} we show an illustration of our method. A shooting method and an Euler discretization are used to simulate the boundary value problem. The blue curve in Figure represents the optimal trajectory interpolating the prescribed points and boundary values.

 The parameters for the trajectory used are $m=1$, $\sigma= 0.5$, $\tau_1 = 1.7$, $\tau_2 = 0.9$, $J=2$, $h=0.01$. Boundary condition are given by: $\gamma(0)=(x(0),y(0),\theta(0))=(2,-1,\pi/3)$, $\gamma(2.28)=(x(2.28),y(2.28),\theta(2.28))=(3.5,5.25,0)$, $v(0)=(a(0), b(0))=(1,1,1)$, $v(2.28)=(a(2.28), b(2.28))$ $=(-1,1,1)$. The interpolation points are $\gamma(1.33)=(1.5,1.5,2\pi/3)$, $\gamma(1.92)=(.5,3.5,\pi/2)$. Figure $3$ (right) shows the evolution of $\theta$ along the time.

In the absence of velocity constraints, the model studied in this example corresponds with a free planar rigid body. The trajectory planning without interpolation points for the obstacle avoidance problem of a planar rigid body was studied using a similar framework  previously by the authors in \cite{BlCaCoCDC} (Section V-A). \end{example}

\section{Final discussion and future research}

We studied the problem of dynamic interpolation for obstacle avoidance on Riemannian manifolds and derived necessary optimality conditions for the trajectory planning problem of mechanical systems specified by a kinetic energy given by a Riemannian metric.

Such optimallity conditions specify a motion of a system along the workspace, interpolating specific points at given times, satisfying boundary conditions, and minimizing an energy functional which depends on an artificial potential function used to avoid static obstacles. Different scenarios were studied: the problem on a Riemannian manifold, the corresponding  sub-Riemannian problem where additional nonholonomic constraints are imposed, systems defined on Lie groups endowed with a left-invariant or bi-invariant Riemannian metrics. Several examples were discussed including left-invariant systems on $SE(3)$, an example on $SO(3)$, and a sub-Riemannian problem on $SE(2)$. All these examples are chosen to cover different aspects of the motion planning problem for several applications in engineering sciences involving Lie group configuration spaces.

The proposed method provides a motion planning algorithm for a class of mechanical control systems that does not require the use of local coordinates in the configuration space.  While we cannot claim rigorously that equation \eqref{3.8} has a solution, given boundary conditions we provide a numerical solution based on Euler's symplectic method which gives a curve that satisfy the necessary and boundary conditions, and  between interpolation points, we solve a boundary value problem by using a shooting method.

The variational approach proposed in this work for the obstacle  avoidance problem allows us  to further study second order optimality conditions for the dynamic interpolation problem and therefore it may be possible to use the approach presented in this work for necessary (first order) conditions to find sufficient (second order) optimality conditions. The existence of global minimizers for the dynamic interpolation problem with obstacle avoidance can be analyzed using similar techniques to the ones developed in \cite{Giambo}, \cite{Giambo2}. 

It is well known that the Pontryagin maximum principle (PMP) can give first order conditions for optimality. As far as we know, such an approach for obstacle avoidance with dynamic interpolation does not exist in the literature. We believe that the study of such a dynamic  interpolation problem from the point of view of PMP, as well as the comparison between both approaches, provides an interesting analysis of the problem discussed in this work.

The study of higher-order variational problems on symmetric spaces and reduction theories for variational problems has attracted considerable interest and has been carried out systematically by several authors. In future work we intend to introduce  interpolation points into such problems and extend the main results presented in this paper to this setting. We will also intend to extend our work to dynamic interpolation for obstacle avoidance with moving obstacles.

\section*{Acknowledgment}

The research of A. Bloch was supported by NSF grant DMS-1613819 and AFOSR. The research of M. Camarinha was partially supported by the Centre for Mathematics of the University of Coimbra -- UID/MAT/00324/2013, funded by the Portuguese Government through FCT/MEC and co-funded by the European Regional Development Fund through the Partnership Agreement PT2020. L. Colombo was supported by MINECO (Spain) grant MTM2016-76072-P.  L. C. and wishes to thank CMUC, Universidade de Coimbra for the hospitality received there where the main part of this work was developed.

\ifCLASSOPTIONcaptionsoff
  \newpage
\fi

\end{document}